\pdfoutput=1
\documentclass[letterpaper, 10 pt, draftcls, onecolumn]{ieeeconf} 
\IEEEoverridecommandlockouts 
\usepackage{graphicx}
\usepackage{subfigure}
\usepackage{amsmath}
\usepackage{amsfonts}
\usepackage[psamsfonts]{amssymb}
\usepackage{cases}
\usepackage{mathtools}
\usepackage{color}
\usepackage{epstopdf}
\usepackage{dsfont}
\usepackage{mathtools}

\newtheorem{definition}{Definition}

\newtheorem{lemma}{Lemma}

\newtheorem{theorem}{Theorem}

\newtheorem{remark}{Remark}

\usepackage{flushend}

\newcommand{\R}{\mathbb{R}}
\newcommand{\N}{\mathbb{N}}

\newcommand{\E}{\mathbb{E}}

\newcommand{\A}{\mathcal{A}}
\newcommand{\aperp}{\A^{\perp}}
\newcommand{\hz}{\hat{z}}
\newcommand{\bone}{\mathds{1}}
\newcommand{\Lhat}{\hat{L}}
\newcommand{\mij}{\Lhat}

\newcommand{\nodeset}{V}
\newcommand{\nodenum}{n}

\mathtoolsset{showonlyrefs=true}

\makeatletter
\renewcommand{\fnum@figure}{Figure \thefigure}
\makeatother

\title{\LARGE \bf
Convergence Rate Estimates for Consensus over Random Graphs
}

\author{Matthew T. Hale and Magnus Egerstedt$^\dag$
\thanks{$^\dag$School of Electrical and Computer Engineering, Georgia Institute of
Technology, Atlanta, GA 30332, USA. Email: \texttt{\{matthale, magnus\}@gatech.edu}.}
}

\begin{document}
\maketitle
\thispagestyle{empty}
\pagestyle{empty}

\begin{abstract}
Multi-agent coordination algorithms with randomized interactions
have seen use in a variety of settings in the multi-agent systems
literature. In some cases, 
these algorithms can be random by design, as in a gossip-like algorithm, and
in other cases they are random due to external factors, as in the case of intermittent communications. 
Targeting both of these scenarios, we present novel convergence rate
estimates for consensus problems solved over random 
graphs. 
Established results provide asymptotic convergence
in this setting, and we provide estimates
of the rate of convergence in two forms. 
First, we estimate decreases in a quadratic Lyapunov
function over time to bound how quickly the agents' disagreement decays, and
second we bound the probability of being at
least a given distance from the point of agreement. 
Both estimates rely on (approximately) computing eigenvalues of the expected matrix exponential
of a random graph's Laplacian, which we do explicitly in terms of the network's
size and edge probability, without assuming that any relationship between them holds.
Simulation results are provided to support the theoretical developments made. 
\end{abstract}

\section{Introduction}
Distributed agreement, often broadly referred
to as the \emph{consensus} problem, is a canonical problem in
distributed coordination and has received attention
in diverse fields such as physics \cite{vicsek95}, signal processing \cite{schizas09},
robotics \cite{ren07},
power systems \cite{nozari14},
and communications \cite{mehyar07}. The goal in such problems is to drive all agents
in a network to a common final state. A key feature of consensus
problems is their distributed nature; consensus is typically carried
out across a network of agents in which each agent communicates with
some other agents, though generally not all of them. The wide range
of fields which study distributed agreement has given rise to corresponding
diversity among consensus problem formulations, and a number of variants
of consensus have been studied in the literature.

In this paper, we derive convergence rates for consensus over
random graphs, studied previously in \cite{hatano05} where asymptotic
convergence was shown. 
In some cases, the motivation for representing a communication network using
random graphs comes from agents using an interaction protocol
that is randomized by design, such as in a gossip-like algorithm \cite{boyd06}.
In other cases, unreliable communications due to poor channel quality,
interference, and other factors can be effectively represented by
a random communication graph \cite{mesbahi10}, and the work
here applies to each of these scenarios. 
This problem formulation has each agent communicating with a random
collection of other agents determined by a random graph. 
Each agent moves toward the average of its neighbors' states,
then a new graph is generated, and each agent moves toward the average state of 
its new neighbors, with this process repeating until convergence is achieved.

We consider networks of a fixed size, and 
we examine consensus over random graphs generated by the Erd\H{o}s-R\'{e}nyi model \cite{erdos59},
in which each possible edge in a graph is present 
with a fixed probability and is 
independent of all other
edges. The Erd\H{o}s-R\'{e}nyi model is used because it accurately captures the 
behavior both of networks with intermittent and unreliable communications \cite{hatano05}
and the behavior of some variants of synchronous gossip algorithms \cite{boyd06}. 
Our approach consists of estimating the expected matrix exponential of 
the Laplacian associated with a random graph and then computing its eigenvalues in
terms of the graph's size and edge probability.
It is shown that the second-largest eigenvalue of this expected matrix exponential
is key in estimating convergence rates, and the main contribution of this paper lies
in explicitly computing this eigenvalue and using it to derive two novel rates of convergence. 
The first estimates the rate of decrease in a quadratic Lyapunov function to bound the rate at which
agents approach agreement in their state values. The second bounds the probability
of the agents' states being at least a given distance away from agreeing on a common state. 
A well-known related work is \cite{touri14}, which develops a general
theory for convergence of consensus over broad classes of time-varying graphs,
in addition to providing general-purpose convergence results.
The emphasis in this paper
is on consensus using the Erd\H{o}s-R\'{e}nyi model due to its prevalence
in some work on multi-agent systems \cite{mesbahi10}. Our work differs from that in \cite{touri14}
because, by fixing the choice 
of graph model, we are able to derive concrete convergence rates and implement them numerically. 

Both estimates rely in some form on computing eigenvalues of random matrices, and
there is an established literature dedicated to doing so \cite{diaconis94,tao12},
including for eigenvalues of random symmetric matrices \cite{alon02,furedi81,wigner58}, and
eigenvalues of random graphs' Laplacians specifically \cite{coja07}.
A common approach to estimating or computing the eigenvalues of a random
symmetric matrix is to let the size of the matrix get arbitrarily large
\cite{diaconis94,furedi81,juhasz82,wigner58}. 
In graph theory, this approach corresponds to letting the number of nodes in a graph
grow arbitrarily large and 
has seen use in spectral graph theory
because it allows one to rigorously state results that hold for almost all graphs \cite{bollobas01}.

In the study of multi-agent systems,
one is often interested in networks of a fixed, small size,
such as in \cite{mehyar07,nozari14,werner06},
and this makes results for asymptotically large networks less applicable in some cases.
As a result, we are motivated to derive convergence results in terms of 
a network's size
without taking it to grow asymptotically large. 
In addition, there are a number of graph theoretic results that
estimate eigenvalues of random graphs' Laplacians when edge probabilities
bear some known relationship to the size of the network, 
e.g., \cite{krivelevich03,tran13}.
In cases where a random
graph is used to model unreliable communications, there is no guarantee that such
a relationship will hold as the quality of communication channels
can depend upon a wide variety of external factors. 
Therefore, we allow edge probabilities to take
values independent of the network's size and state our results in terms
of both a network's size and its edge probability without making assumptions about either. 
Our results also do not require the common assumption that unions of agents' communication
graphs are connected across finite
intervals of a prescribed length \cite{jadbabaie03}, which can be difficult to enforce and verify
in some practical settings.

The rest of the paper is organized as follows. Section II reviews the necessary
background on consensus and graph theory, including that on random graphs.
Next, Section III establishes asymptotic convergence of consensus over random graphs.
Then Section IV presents our main
results on the rate of convergence of consensus over random graphs. 
Section V then presents simulation
results to verify the convergence rates we present. Finally, Section VI concludes the paper.

\section{Review of Graph Theory and Consensus}
In this section, we review the basic elements of graph
theory required for the remainder of the paper.
We introduce unweighted, undirected graphs, then review
the consensus problem, and finally review random graphs. 

\subsection{Basic Graph Theory}
All graphs in this paper are assumed to be unweighted and
undirected. 
Such graphs are defined by pairwise relationships over a finite set
of nodes or vertices. Suppose that a graph has a set $\nodeset$ 
of $\nodenum$ vertices,
with $\nodenum \in \N$, and index these vertices over the set
$\{1, \ldots, \nodenum\}$. We define the \emph{edge set}
\begin{equation}
E \subset V \times V,
\end{equation}
and say there is an edge between nodes $i$ and $j$
if $(i, j) \in E$. 
A graph $G$ is then formally defined 
as the $2$-tuple $G = (V, E)$. 
Throughout this paper, all edges are undirected and
an edge $(i, j) \in E$ is not distinguished from
the edge $(j, i) \in E$. We do not allow self loops and
therefore $(i, i) \not\in E$ for all $i$ and all
graphs $G$. 

The degree associated with node $i$ is defined as the total
number of edges that connect node $i$ to some other node.
Using $|\cdot|$ to denote the cardinality of a set
and using $d_i$ to denote the degree of node $i$, we have
\begin{equation}
d_i = \big|\{j \mid (i, j) \in E\}\big|.
\end{equation}
The $n \times n$ degree matrix associated with a graph $G$ is then defined
as
\begin{equation}
  D(G) = \left(\begin{array}{cccc} d_1 & 0      & \cdots & 0 \\
                                     0 & d_2    & \cdots & 0 \\
                               \vdots  & \vdots & \ddots & \vdots \\
                                    0  &      0 & \cdots & d_n \end{array}\right),
\end{equation}
which will be written simply as $D$ when the graph $G$ is understood.

The $n \times n$ adjacency matrix associated with
$G$, denoted $A(G)$, is defined element-wise as
\begin{equation}
  a_{i,j} = \begin{cases} 1 & (i, j) \in E \\
                         0 & \textnormal{otherwise} \end{cases}
\end{equation}
where $a_{i,j}$ is the $(i,j)$th entry in $A(G)$. 
Note that, because $(i, j) \in E$ implies $(j, i) \in E$, $A(G)$
is a symmetric matrix. In addition, the absence of self loops
results in $A(G)$ having zeroes on its main diagonal for all
graphs $G$. We will simply write $A$ when $G$ is clear from
context. 


Using $D(G)$ and $A(G)$, we define the Laplacian associated with a graph $G$ as
\begin{equation}
L(G) = D(G) - A(G),
\end{equation}
and note that 
the Laplacian of any undirected, unweighted graph is a symmetric, positive semi-definite
matrix \cite{godsil01}. In particular, the Laplacians of such graphs
have all non-negative eigenvalues.
Letting $\lambda_i(M)$ denote the $i^{th}$ smallest eigenvalue
of a matrix $M$, for any graph $G$ it is known that $\lambda_1\big(L(G)\big) = 0$ and
that $\bone = (1, 1, \ldots, 1)^T$ is an eigenvector associated with $\lambda_1$ \cite{mesbahi10},
i.e., that $\bone$ is in the nullspace of the Laplacian of any undirected, unweighted graph.  
For any graph Laplacian $L$ of size $n \times n$, we find
\begin{equation} \label{eq:eiglist}
0 = \lambda_1 \leq \lambda_2 \leq \cdots \leq \lambda_n.
\end{equation}

A graph $G$ is said to be \emph{connected} if, for any two nodes $i$ and $j$ in the graph,
there exists a sequence of edges one can traverse to travel from node $i$ to node $j$. 
A seminal result in graph theory provides that $G$ is connected if and only
if $\lambda_2 > 0$ \cite{fiedler73}.

\subsection{The Consensus Protocol} \label{ss:reviewcon}
A canonical problem in multi-agent control is that of consensus. Consensus
consists of having a collection of agents, e.g., robots or
mobile sensor nodes, agree on a common value in a distributed way.
The term distributed refers to the fact that each agent in a network
can only communicate with some other agents in the network, but 
in general not all other agents.  
Suppose each agent has a scalar state\footnote{All results in this paper
are easily extended to the case of agents having vectors of states
by replacing $L$ in the consensus protocol in Equation~\eqref{eq:con1} with $I \otimes L$, where
$\otimes$ denotes the Kronecker product of two matrices. 
We focus on the scalar case to simplify issues of dimensionality.}, with agent $i$'s state denoted
$x_i$, and assemble these into the ensemble state vector
\begin{equation}
x = \left(\begin{array}{c} x_1 \\ \vdots \\ x_n \end{array}\right)^T \in \R^n. 
\end{equation}

If one represents the
agents' communications using a graph $G$ (where edges connect those agents that
exchange information), then the continuous-time consensus protocol \cite{olfati-saber04} takes the form
\begin{equation} \label{eq:con1}
\dot{x} = -Lx,
\end{equation}
where $L$ is the Laplacian of the graph $G$. 

The following well-studied theorem from \cite{olfati-s03} establishes 
convergence of the consensus protocol in continuous time when $G$ is connected. 
In it, the symbol
\begin{equation}
\bar{x}(0) := \frac{1}{n}\sum_{i=1}^{n} x_i(0)
\end{equation}
is used to denote the centroid of the agents' initial
states (which is a scalar) and $\bar{x}(0)\bone$ 
is used to denote the vector
in $\R^n$ whose entries are all equal to $\bar{x}(0)$. 
\begin{theorem} \label{thm:consensus}
Let $G$ be a connected graph and let an initial ensemble state $x(0)$ be given. Then the consensus protocol
\begin{equation} \label{eq:consensus}
\dot{x} = -Lx
\end{equation}
asymptotically converges element-wise to the centroid of the agents' initial states, 
i.e., $x(t) \to \bar{x}(0)\bone$. In addition, the rate at which $x(t)$
approaches $\bar{x}(0)\bone$ is governed by $\lambda_2$. 
\end{theorem}
\emph{Proof}: See \cite{olfati-s03}. \hfill $\blacksquare$

Because the goal in consensus is for a team of nodes to reach a common value, one
often refers to the agreement subspace of the agents, defined below as
was done in \cite{hatano05}. 

\begin{definition} (\cite{hatano05}) \label{def:agree}
The agreement subspace is defined as the set of all points at which all agents have
the same state value, i.e., $x_i = x_j$ for all $i$ and $j$. Formally, it is defined as 
\begin{equation}
\mathcal{A} = \textnormal{span}\{\mathds{1}\},
\end{equation}
where $\mathds{1}$ is the vector of all ones in $\R^n$. \hfill $\triangle$
\end{definition}

In Section~\ref{sec:conv}, 
we will show asymptotic convergence to $\mathcal{A}$ for consensus over random
graphs by showing that the agents' disagreement goes to zero. 
Next, we introduce the random graph model used in the remainder of the paper. 

\subsection{Random Graphs}
A common model for random graphs is the Erd\H{o}s-R\'{e}nyi model,
originally published in \cite{erdos59}, and we use it here because it
accurately captures the behavior of two cases of interest. First,
some algorithms are randomized by design, such as gossip algorithms \cite{boyd06},
and Erd\H{o}s-R\'{e}nyi graphs can model the behavior of such algorithms
in some cases. Second, members of a network sometimes share information over
communication channels which are intermittently lost and regained, 
and this behavior is well-modeled by Erd\H{o}s-R\'{e}nyi graphs \cite{hatano05}. 
This model takes two parameters
to generate random graphs: a number of nodes $n \in \N$ and
an edge probability\footnote{The cases $p = 0$ and
$p = 1$ provide edgeless graphs and complete graphs, respectively,
and are omitted because their behavior is deterministic.} $p \in (0, 1)$.
The Erd\H{o}s-R\'{e}nyi model generates graphs on $n$ nodes
whose edge sets contain
each possible edge with probability $p$, independent of all other
edges. 
Formally, for each admissible
$i$ and $j$, we have
\begin{equation}
\mathbb{P}[(i, j) \in E] = p.
\end{equation}
An alternative characterization that we use later can be stated
in terms of the elements of the adjacency matrix of a random graph:
for $n$ nodes and edge probability $p$, we find
\begin{equation}
\mathbb{P}[a_{i,j} = 1] = \mathbb{P}[a_{j,i} = 1] = p \textnormal{ and } \mathbb{P}[a_{i,j} = 0] = \mathbb{P}[a_{j,i} = 0] = 1 - p.
\end{equation}
Thus each $a_{i,j}$ is a Bernoulli random variable. 

We use $\mathcal{G}(n, p)$ to denote the sample space of all possible random graphs
generated by the Erd\H{o}s-R\'{e}nyi model on $n$ nodes with edge probability $p$, and
we use $\mathcal{L}(n, p)$ to denote the set of Laplacians of all such graphs. 
One approach to analyzing random graphs that has seen use in the graph theory
literature consists of taking the limit as $n$ goes to infinity, 
with the benefit of this approach being the ability to 
rigorously determine which properties hold for almost all graphs.
However, motivated by the study of multi-agent systems, we are interested in
networks of fixed size and therefore develop our results 
in terms of a fixed value of $n$. 
Toward doing so, we show asymptotic convergence of consensus over random
graphs in the next section. 


\section{Consensus over Random Graphs} \label{sec:conv}
In this section we examine consensus where the agents' communication
graph at each timestep is a random graph. 
We then show asymptotic convergence of this update law in order to help establish
the role of convergence rates in its analysis. 
The theoretical results of this section are not new, but are presented to
contextualize the remainder of the paper. 
This section closely follows the approach of \cite{hatano05}
where these results were originally published. 

\subsection{Time-Sampled Consensus} \label{ss:3a}
We assume that all communication graphs
hold constant for some positive amount of time $\delta > 0$  and we seek to
examine the evolution of the system's state under this condition.
We follow the problem setup of \cite{hatano05} in which states generated by
the consensus protocol in Equation~\eqref{eq:consensus} are sampled in time by defining
the state $z(k)$ as
\begin{equation}
z(k) = x(k\delta),
\end{equation} 
where the communication graph among the agents is assumed to hold
constant over the interval $[k\delta, (k+1)\delta)$.
We note that this problem is distinct from consensus performed in discrete time;
the problem we consider 
analyzes samples of the continuous-time state $x(t)$ rather than having
states actually evolve in discrete time as in \cite{jadbabaie03}. 
All agents still execute the protocol $\dot{x} = -Lx$ and,
as in \cite{hatano05}, $z(k)$ is used only as a theoretical
tool for analyzing the behavior of the continuous state $x(t)$ over time
as the agents' communication graphs change. It is because these communication
graphs hold constant over $[k\delta, (k+1)\delta)$ that we analyze
the agents' states only at these points in time. 

At each timestep $k$, the system will generally have a different graph
Laplacian than it had at time $k-1$. We denote the communication graph active
at time $k$ by $G_k$ and denote its Laplacian by $L_k$. 
The solution to Equation~\eqref{eq:consensus} is $x(t) = e^{-Lt}x(0)$,
and for any times $t_1$ and $t_2$ with $t_1 < t_2$ we find
\begin{equation}
x(t_2) = e^{-L(t_2 - t_1)}x(t_1). 
\end{equation}
Setting $t_2 = (k+1)\delta$ and $t_1 = k\delta$ then gives
\begin{equation} \label{eq:xddef}
x\big((k+1)\delta\big) = e^{-\delta L_k}x(k\delta)
\end{equation}
because the graph $G_k$ is constant over the interval
$[k\delta, (k+1)\delta)$. Equation \eqref{eq:xddef} itself is equal to
\begin{equation} \label{eq:discon}
z(k+1) = e^{-\delta L_k}z(k)
\end{equation}
and this is the protocol that we analyze in the remainder of the paper. 

We emphasize that the agents' states still evolve in continuous time
and that we analyze the continuous time signal $x(t)$ by analyzing
samples taken every $\delta$ seconds. In Section \ref{sec:rates} we 
provide convergence rates for Equation \eqref{eq:discon}.
Toward doing so, we next show that Equation~\eqref{eq:discon} converges 
asymptotically to $\mathcal{A}$.  

\subsection{Establishing Asymptotic Convergence} \label{ss:asymptotic}
We will assess convergence of consensus over random graphs by showing 
convergence to the agreement set $\mathcal{A}$, defined in Definition~\ref{def:agree}. Define the
distance from a point $y \in \R^n$ to $\mathcal{A}$ as
\begin{equation} \label{eq:distdef}
\textnormal{dist}(y, \mathcal{A}) = \inf_{z \in \mathcal{A}} \|z - y\|_2. 
\end{equation}
As in \cite{hatano05}, define
the orthogonal complement of the agreement subspace as
\begin{equation}
\aperp = \{x \mid x^Ta = 0 \textnormal{ for all } a \in \A\}.
\end{equation}
In particular $x^T\bone = 0$ for all $x \in \aperp$.
Then define the Euclidean projection of $z(k)$ onto $\aperp$ as
\begin{equation}
\hat{z}(k) = \Pi_{\aperp}[z(k)] = \left(I - \frac{1}{n}J\right)z(k),
\end{equation}
where $J$ is the $n \times n$ matrix of ones, and 
where $\hat{z}(k)$ captures the disagreement among agents by excluding
the part of $z(k)$ that lies in $\mathcal{A}$. 

It was noted in Theorem~\ref{thm:consensus} that the consensus protocol 
over static graphs converges
element-wise to centroid of the agents' initial states. 
The consensus protocol in Equation \eqref{eq:discon} also converges to the centroid
of the agents' initial states, despite being run over random graphs. 
To see why, first observe that the consensus protocol in Equation \eqref{eq:discon}
has converged when $z(k+1) - z(k) = 0$, i.e., when
\begin{equation}
e^{-\delta L_k}z(k) - z(k) = (e^{-\delta L_k} - I)z(k) = 0 \label{eq:usenull}. 
\end{equation}
Because $\mathcal{A}$ is a subspace, we can decompose $z(k)$ into
two parts according to 
\begin{equation}
z(k) = \Pi_{\mathcal{A}}[z(k)] + \Pi_{\aperp}[z(k)] = \Pi_{\mathcal{A}}[z(k)] + \hz(k). 
\end{equation}
Substituting this decomposition into Equation \eqref{eq:usenull} we find
\begin{align}
 (e^{-\delta L_k} - I)z(k) &= \Big(-\delta L_k + \frac{1}{2}\delta^2L_k^2 - \frac{1}{3!}\delta^3L_k^3 + \cdots\Big)\big(\Pi_{\mathcal{A}}[z(k)] + \hz(k)\big) \\
                           &= \left(-\delta L_k + \frac{1}{2}\delta^2L_k^2 - \frac{1}{3!}\delta^3L_k^3 + \cdots\right)\hz(k) \\
                           &= 0 \label{eq:usenull2},
\end{align}
where we have used the fact that $\Pi_{\mathcal{A}}[z(k)]$ is in the nullspace of  
all $L_k \in \mathcal{L}(n, p)$ by definition. 
For Equation \eqref{eq:discon} to have converged, one must therefore have $\hz(k)$ in the nullspace
of all $L_k \in \mathcal{L}(n, p)$. Since $\hz(k) \perp \mathcal{A}$, we find $\hz(k) = 0$.
Then if Equation \eqref{eq:discon} has converged at time $k$, it has converged to $\Pi_{\mathcal{A}}[z(k)]$. 

To see that $\Pi_{\mathcal{A}}[z(k)] = \bar{z}(0)$, consider $\Pi_{\mathcal{A}}[z(1)]$.
We have
\begin{equation}
\Pi_{\mathcal{A}}[z(1)] = z(1) - \hz(1) = z(1) - \left(I - \frac{1}{n}J\right)z(1) = \bar{z}(1)\bone. 
\end{equation}
Examining $\bar{z}(1)$, we see that
\begin{align}
\bar{z}(1) &= \frac{1}{n}\bone^Tz(1) = \frac{1}{n}\bone^Te^{-\delta L_0}z(0) \\
           &= \frac{1}{n}\bone^T\left(I - \delta L_0 + \frac{\delta^2 L_0^2}{2} - \cdots\right)z(0) \\
           &= \frac{1}{n}\bone^Tz(0) = \bar{z}(0),
\end{align}
where we have used $\bone^TL_0 = (L_0^T\bone)^T = (L_0\bone)^T = 0$ because $L_0$ is
symmetric and $\bone$ is in the nullspace of $L_0$. A simple inductive argument shows
that $\bar{z}(k) = \bar{z}(0)$ for all $k$. 

Therefore, the distance from $z(k)$ to $\mathcal{A}$
is equal to that from $z(k)$ to $\bar{z}(0)\bone$. 
Noting that $\bar{z}(0)\bone$, 
whose entries are all $\bar{z}(0)$, 
is equal to $\frac{1}{n}\mathds{1}^Tz(0)\mathds{1}$, we have
\begin{align}
\textnormal{dist}\big(z(k), \mathcal{A}\big)^2 &= \left\|z(k) - \frac{1}{n}\mathds{1}^Tz(0)\mathds{1}\right\|^2 \\
                                               &= z(k)^Tz(k) - n\bar{z}(0)^2 \\
                                               &= \frac{1}{n}z(k)^T\Lhat z(k), \label{eq:vnice}
\end{align}
where we have used $\Lhat := nI - J$, with $I$ the $n \times n$ identity matrix and $J$ 
the $n \times n$ matrix of ones. 

In light of the form of $\textnormal{dist}\big(z(k), \mathcal{A}\big)$, we will 
show convergence of the protocol in Equation~\eqref{eq:discon}
using the quadratic Lyapunov function
\begin{equation} \label{eq:Vdef}
V\big(z(k)\big) = \frac{1}{n}z(k)^T\Lhat z(k).
\end{equation}

We also state the following definition which will be used to characterize stochastic convergence.
\begin{definition} \label{def:wp1}
A random sequence $\{y(k)\}$ in $\R^n$ converges to $y \in \R^n$ with
probability $1$ if,
for every $\epsilon > 0$,
\begin{equation}
\mathbb{P}\left[\sup_{N \leq k < \infty} \|y(k) - y\|_2 \geq \epsilon\right] \to 0 \textnormal{ as } N \to \infty.
\end{equation}
\hfill $\triangle$
\end{definition}

We now have the following theorem that proves asymptotic convergence of
the consensus protocol in Equation~\eqref{eq:discon}, due originally
to \cite{hatano05}. 

\begin{theorem} \label{thm:asymptotic}
Fix a number of nodes $n \in \N$ and an edge probability $p \in (0, 1)$, and let
$z(0)$ be given. 
The consensus protocol 
\begin{equation} \label{eq:thmcon}
z(k+1) = e^{-\delta L_k}z(k),
\end{equation}
where $L_k \in \mathcal{L}(n, p)$
for all $k$, converges to $\mathcal{A}$ with probability $1$. 
In addition,
\begin{equation}
\hz(k)^T\E[e^{-2\delta L_k} - I]\hz(k) \to 0
\end{equation}
with probability $1$ and 
\begin{equation}
\mathbb{P}\left[\sup_{N \leq k < \infty} \|\hz(k)\|_2^2 \geq \gamma\right] \leq \frac{\hz(0)^T\hz(0)}{\gamma}\lambda_{n-1}\big(\E[e^{-2\delta L_k}]\big)^N
\end{equation}
where $\lambda_{n-1}(\cdot)$ denotes the second-largest eigenvalue of a matrix. 
\end{theorem}
\emph{Proof:} See \cite{hatano05}. 
\hfill $\blacksquare$

More details on this theorem can be found in \cite{hatano05}, 
where its proof originally appeared, and in \cite{mesbahi10}. 
Having established asymptotic convergence, we derive rates of convergence for consensus over random graphs in the
next section.

\section{Convergence Rates for Consensus over Random Graphs} \label{sec:rates}
Section~\ref{sec:conv} showed that
consensus over random graphs converges asymptotically, and 
in this section we derive our main results on two rates of convergence
for consensus over random graphs. Section \ref{ss:crates} is based on
work in \cite{hatano05}, though the form of convergence rate we derive
is different from the one derived there. The remainder of the section
then presents our main results on novel, explicit convergence rates. 

\subsection{Convergence Rates} \label{ss:crates}
We now highlight the utility of convergence rates in the context of Theorem 
\ref{thm:asymptotic}. 
To do so, we examine the evolution of the disagreement among agents. 
Using the same Lyapunov function as in Equation \eqref{eq:Vdef},
we find
\begin{align}
\E[V\big(\hz(k+1)\big) - V\big(\hz(k)\big) \mid \hz(k)] &= \hz(k)^T\E[e^{-2\delta L_k} - I]\hz(k) \\
                                                        &= \hz(k)^T\E[e^{-2\delta L_k}]\hz(k) - \hz(k)^T\hz(k).
\end{align}
By definition, $\hz(k)$ is orthogonal to $\bone$.  
$\bone$ is also the eigenvector
associated with the largest eigenvalue of $e^{-2\delta L_k}$
for all $L_k \in \mathcal{L}(n, p)$ 
(because the eigenvalues of 
$e^{-2\delta L_k}$ are $e^{-2\delta\lambda_i}$ for
each $\lambda_i$ in Equation \eqref{eq:eiglist}). 
Then  we find that
\begin{equation}
\hz(k)^T\E\big[e^{-2\delta L_k}\big]\hz(k) \leq \lambda_{n-1}\left(\E\big[e^{-2\delta L_k}\big]\right)\|\hz(k)\|_2^2,
\end{equation}
where $\lambda_{n-1}(M)$ denotes the second largest eigenvalue of an $n \times n$ matrix $M$. 
Consequently, we have
\begin{equation} \label{eq:lamn-1}
\E[V\big(\hz(k+1)\big) - V\big(\hz(k)\big) \mid \hz(k)] \leq \big(\lambda_{n-1}(\E[e^{-2\delta L_k}]) - 1\big)\|\hz(k)\|_2^2,
\end{equation}
and by Theorem~\ref{thm:asymptotic} we have
\begin{equation} \label{eq:probout}
\mathbb{P}\left[\sup_{N \leq k < \infty} \|\hz(k)\|_2^2 \geq \gamma\right] \leq \frac{\hz(0)^T\hz(0)}{\gamma}\lambda_{n-1}\big(\E[e^{-2\delta L_k}]\big)^N.
\end{equation}

Therefore, both convergence rates depend upon $\lambda_{n-1}(\E[e^{-2\delta L_k}])$,
and we compute it next. 

\subsection{Computing $\lambda_{n-1}\big(\E[e^{-2\delta L_k}]\big)$}
We have the following lemma concerning eigenvalues of a matrix of the form $a I + b (J - I)$, which
we will use below.
\begin{lemma} \label{lem:ab}
A matrix $M$ of the form $a I + b (J - I)$, namely
\begin{equation}
M = \left(\begin{array}{cccc} a    &      b & \cdots &    b \\
                              b    &      a & \cdots &    b \\
                            \vdots & \vdots & \ddots & \vdots \\
                              b    &      b & \cdots &    a \end{array}\right) \in \R^{n \times n},
\end{equation}
 has $a + (n-1)b$
as an eigenvalue with multiplicity one and $a - b$ as an eigenvalue with multiplicity $n - 1$.
\end{lemma}
\emph{Proof:} We proceed using a series of row operations that will preserve the 
characteristic polynomial of $M$. We see that
\begin{equation}
|M - \lambda I| = \left|\left(\begin{array}{cccc} a - \lambda     &               b  &      \cdots &    b \\
                                                             b    &      a - \lambda &      \cdots &    b \\
                                                 \vdots           & \vdots           &      \ddots & \vdots \\
                                                   b              &      b           &      \cdots &    a - \lambda \end{array}\right)\right|.
\end{equation}
Next, we add rows $2$ through $n$ to row $1$, giving
\begin{equation}
|M - \lambda I| = \big(a + (n-1)b - \lambda\big)\left|\left(\begin{array}{cccc} 1   &            1     &     \cdots  &    1 \\
                                                                               b    &      a - \lambda &     \cdots  &    b \\
                                                                             \vdots &          \vdots  &     \ddots  & \vdots \\
                                                                               b    &      b           &     \cdots  &    a - \lambda \end{array}\right)\right|.
\end{equation}

Subtracting $b$ times row $1$ from each other row, we find
\begin{equation}
|M - \lambda I| = \big(a + (n-1)b - \lambda\big)\left|\left(\begin{array}{cccc} 1   &   1                  &     \cdots  &    1 \\
                                                                               0    &      a - b - \lambda &     \cdots  &    0 \\
                                                                             \vdots & \vdots               &     \ddots  & \vdots \\
                                                                               0    &      0               &     \cdots  &    a - b - \lambda \end{array}\right)\right|,
\end{equation}
where the matrix on the right-hand side is upper-triangular. The determinant on the right-hand side is then the product of the
diagonal entries of that matrix, resulting in
\begin{equation}
|M - \lambda I| = (a + (n-1)b - \lambda)(a - b - \lambda)^{n-1},
\end{equation}
whose roots are indeed $a + (n-1)b$, with multiplicity $1$, and $a - b$ with multiplicity $n-1$.  \hfill $\blacksquare$

We now derive the expected value of the first four powers of a random graph's Laplacian, and
below we will use these results to approximate the Laplacian's expected matrix exponential. 
As above, we use the notation $\Lhat = nI - J$. 

\begin{lemma} \label{lem:lpow}
Let a number of nodes $n \in \N$ and edge probability $p \in (0, 1)$ be given. 
The Laplacian $L \in \mathcal{L}(n, p)$ of a graph $G \in \mathcal{G}(n, p)$
satisfies
\begin{align}
\E[L]   &= p\mij \\
\E[L^2] &= \left[(n-2)p^2 + 2p\right]\mij \\
\E[L^3] &= \left[(n-2)(n-4)p^3 + 6(n-2)p^2 + 4p\right]\mij \\
\E[L^4] &= \big[(n-7)(n-3)(n-2)p^4 + 6(2n-7)(n-2)p^3 + 25(n-2)p^2 + 8p\big]\mij.
\end{align}
\end{lemma}
\emph{Proof sketch:} 
We sketch the proof to avoid exposition on many tedious computations
and instead elaborate on the core arguments used to derive the above results.

The general form of graph Laplacian for $G \in \mathcal{G}(n, p)$ is
\begin{equation}
\left(\begin{array}{cccc} 
          \sum_{\substack{j=1 \\ j \neq 1}}^{n}a_{1, j}             & -a_{1, 2}                                                 & \cdots & -a_{1, n} \\
          -a_{1, 2}                                                 & \sum_{\substack{j=1 \\ j \neq 2}}^{n}a_{2, j}             & \cdots & -a_{2, n} \\
           \vdots                                                   &   \vdots                                                  & \ddots & \vdots    \\
          -a_{1, n}                                                 & -a_{2, n}                                                 & \cdots & \sum_{\substack{j=1 \\ j \neq n}}^{n}a_{n, j} \end{array}\right),
\end{equation}
where each term $a_{i, j}$ is a Bernoulli random variable
with expectation equal to $p$. 
The off-diagonal entries of $L$ have $\E[L_{ij}] = \E[-a_{i,j}] = -p$, while
linearity of $\E[\cdot]$ gives $\E[L_{ii}] = (n-1)p$ for diagonal
entries of $L$. From this we find $\E[L] = (n-1)pI - p(J - I) = p\Lhat$. 

Computing the expectation of the entries of $L^2$ requires one to consider two cases. The diagonal entries
of $L^2$ are formed by the product of row $i$ of $L$ with column $i$ of $L$ (and by symmetry of $L$ these are identical),
while the off-diagonal entries result from the product of row $i$ and column $j$ of $L$ (which are not identical
when $i \neq j$). 
It is important to note that $a_{i, j}^2 = a_{i, j}$ because $a_{i, j}$ is a Bernoulli random
variable. As a result, when computing expectations one finds that $\E[a_{i,j}^k] = p$ for all $k \in \N$, while
products of $\ell$ distinct off-diagonal entries of $A$ have expectation equal to $p^\ell$. 
In the case of $\E[L^2]$, a diagonal entry takes the form
\begin{align}
\E\left[(L^2)_{ii}\right] &= \E\left[\sum_{\substack{p = 1 \\ p \neq i}}^{n} a_{i,p}^2 + \left(\sum_{\substack{q = 1 \\ q \neq i}}^{n} a_{i,q}\right)^2\right] \\
                          &= (n-1)(n-2)p^2 + 2(n-1)p,
\end{align}
while an off-diagonal entry takes the form
\begin{align}
\E[(L^2)_{ij}] &= \E\left[\sum_{\substack{k=1 \\ k \neq i, j}}^{n} a_{i,k}a_{k,j}\right] - \E\left[a_{i,j}\sum_{\substack{p = 1 \\ p \neq i}}^{n} a_{i,p}\right] - \E\left[a_{i,j}\sum_{\substack{q = 1 \\ q \neq i}}^{n} a_{q,j}\right] \\
               &= -(n-2)p^2 - 2p.
\end{align}
Then we find
\begin{align}
\E[L^2] &= \big((n-1)(n-2)p^2 + 2(n-1)p\big)I + \big({-(n-2)}p^2 - 2p\big)(J - I) \\
        &= \big((n-2)p^2 + 2p\big)\Lhat,
\end{align}
as above. 

Computing the general form of $L^3$ is done by multiplying the general form of $L^2$ by that of $L$ and
the general form of $L^4$ is found by squaring the general form of $L^2$. Having found these two general forms,
one follows the above strategy for computing their expected values: first replace $a_{i, j}^k$ with
$a_{i, j}$ for all $k > 1$ and then compute the expectation of the products of $\ell$ distinct
off-diagonal entries of $A$ as $p^{\ell}$, resulting in the above. \hfill $\blacksquare$

We use the first four powers of $L$ because that is all that is required for accurate convergence rate
estimates, as will be shown in Section~\ref{sec:simulation}.
One way in which this accuracy is attained is through the choice of $\delta$. Many choices of $\delta$
are possible and we choose $\delta = 1/n$, for several reasons. First, it was assumed in Section
\ref{ss:3a} that the communication graph of the system is constant between samples, i.e., 
that $G_k$ does not change over the interval $[k\delta, (k+1)\delta)$. 
As a network grows,
the number of possible edges does too and thus a larger network has more ways in which its communication
topology may change. As a result, a larger network should use shorter sampling times and
$\delta$ should decrease as $n$ increases. The choice of $\delta = 1/n$
provides a simple means of enforcing this condition. 

Second, the use of $\delta = 1/n$ is also partly inspired
by the same choice made in \cite{jadbabaie03} for discrete-time consensus where it is a necessary condition
for stability; though not necessary for stability in the current paper, we make the same choice 
 to help retain the same broad applicability of the results in \cite{jadbabaie03}. 
Using this choice of $\delta$, we now give the approximate value of $\lambda_{n-1}(\E[e^{-2\delta L_k}])$
in terms of $n$ and $p$. 

\begin{theorem} \label{thm:lambda}
Suppose $\delta = 1/n$ and define
\begin{align}
\kappa_1 &:= p \\
\kappa_2 &:= (n-2)p^2 + 2p \\
\kappa_3 &:= (n-2)(n-4)p^3 + 6(n-2)p^2 + 4p \\
\kappa_4 &:= (n-7)(n-3)(n-2)p^4 + 6(2n-7)(n-2)p^3 + 25(n-2)p^2 + 8p \\
\mu(n, p) &:= -2\frac{\kappa_1}{n} + 2\frac{\kappa_2}{n^2} - \frac{4}{3}\frac{\kappa_3}{n^3} + \frac{2}{3}\frac{\kappa_4}{n^4}.
\end{align}
Then
\begin{equation}
E[e^{-2\delta L_k}] \approx I + \mu(n, p)\mij,
\end{equation}
and consequently we have
\begin{equation}
\lambda_{n-1}(\E[e^{-2\delta L_k}]) \approx 1 + n\mu(n, p).
\end{equation}
\end{theorem}
\emph{Proof:} Taylor expanding the matrix exponential, we find
\begin{equation} \label{eq:maineq1}
\E[e^{-2\delta L_k}] \approx I - 2\delta\E[L_k] + 2\delta^2\E[L_k^2] - \frac{4}{3}\delta^3\E[L_k^3] + \frac{2}{3}\delta^4\E[L_k^4],
\end{equation}
where the exact matrix exponential is well-approximated by this truncation in part due
to the aforementioned choice of $\delta = 1/n$. Substituting $\delta = 1/n$ and the
results of Lemma~\ref{lem:lpow} into Equation~\eqref{eq:maineq1} gives
\begin{align} \label{eq:maineq2}
\E[e^{-2\delta L_k}] &\approx I + \left(-2\frac{\kappa_1}{n} + 2\frac{\kappa_2}{n^2} - \frac{4}{3}\frac{\kappa_3}{n^3} + \frac{2}{3}\frac{\kappa_4}{n^4}\right)\Lhat \\
                     &= I + \mu(n, p)\Lhat,
\end{align}
which establishes the first part of the theorem. 

Next, we use the definition of $\Lhat$ as $nI - J$ to find
\begin{equation}
I + \mu(n, p) \Lhat = \big(1 + (n-1)\mu(n, p)\big)I - \mu(n, p)\big(J - I\big). 
\end{equation}
Using Lemma~\ref{lem:ab} with $a = 1 + (n-1)\mu(n, p)$ and $b = -\mu(n, p)$, we find
that the largest eigenvalue of $\E[e^{-2\delta L_k}]$ is $1$, 
and the second largest through smallest eigenvalues of $\E[e^{-2\delta L_k}]$ are all approximately
$1 + n\mu(n, p)$. In particular, $\lambda_{n-1}\big(\E[e^{-2\delta L_k}]\big) \approx 1 + n\mu(n, p)$,
as desired. 
\hfill $\blacksquare$

\subsection{Explicit Convergence Rates for Consensus over Random Graphs}
We now present our unified main convergence rates for consensus over random graphs,
stated in terms of network size $n$ and edge probability $p$. 

\begin{figure}
\centering
\includegraphics[width=3.6in]{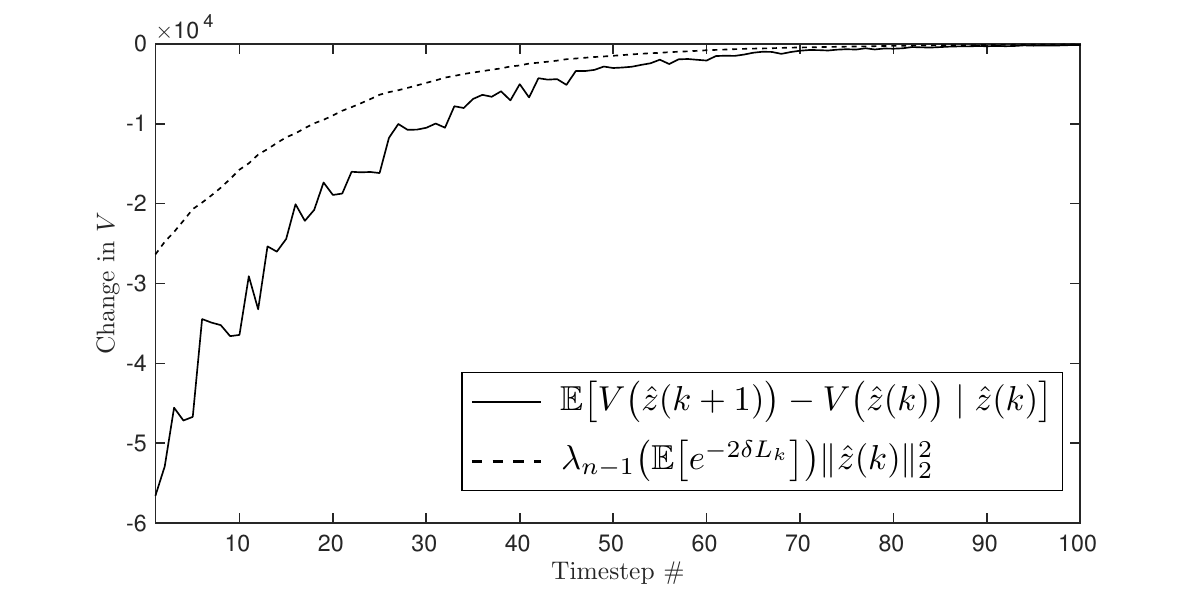}
\caption{A simulation run showing the expected decrease in $V$ (solid line) 
and its theoretical upper bound (dashed line) as in Equation~\protect\eqref{eq:sim1}. Here we see that the upper
bound indeed holds across the whole time horizon, with the dashed line always above
the solid line. In addition, the upper bound
becomes more accurate across iterations, with it accurately predicting when
the expected decrease in $V$ is near zero. 
}
\label{fig:vdiff}
\end{figure}

\begin{theorem} \label{thm:main}
Let a network size $n \in \N$ and an edge probability $p \in (0, 1)$ be given, and
let $\mu(n, p)$ be as defined in Theorem~\ref{thm:lambda}.
For sampling constant $\delta = 1/n$ and $\gamma > 0$ we have
\begin{equation} \label{eq:mainprob}
\mathbb{P}\left[\sup_{N \leq k < \infty} \|\hz(k)\|_2^2 \geq \gamma\right] \leq \frac{\hz(0)^T\hz(0)}{\gamma}\big(1 + n\mu(n, p)\big)^N
\end{equation}
for all $N \in \N$.
In addition, the expected decrease in the
Lyapunov function $V(z) = \frac{1}{n}z^T\Lhat z$ from time $k$ to time $k+1$ is bounded according to
\vspace{-0.3cm}
\begin{equation} \label{eq:nmu}
\E\big[V\big(\hz(k+1)\big) - V\big(\hz(k)\big) \mid \hz(k)\big] \leq n\mu(n, p)\|\hz(k)\|_2^2.
\end{equation}
\end{theorem}
\emph{Proof:} This follows from Equations~\eqref{eq:lamn-1} and \eqref{eq:probout} and
Theorem~\ref{thm:lambda}. 
\hfill $\blacksquare$

\begin{remark}
In truncating the Taylor expansion of $\E[e^{-2\delta L_k}]$, some level of error
in inevitably introduced into the value of $\lambda_{n-1}\big(\E[e^{-2\delta L_k}]\big)$
that results. By truncating after the fourth term, the error in this case is on
the order of $\frac{2^5}{5!}\delta^5\E[L_k^5]$. For the choice of $\delta = 1/n$,
this gives error on the order of $2^5/5!$. This error can be further mitigated
by other choices of $\delta$, though the choice of $1/n$ will suffice in many cases.
In general, the approximations we make are more accurate for smaller values of $p$
because such values cause higher-order terms in the expansion of
$\E[e^{-2\delta L_k}]$ to be dominated by the lower-order terms we include
in our approximations. \hfill $\triangle$
\end{remark}

The appeal of using these convergence rate estimates together is that one need only
compute the constant $\mu(n, p)$ and then
both estimates can be used. Both provide information at each timestep because
they rely on the current iteration count and can therefore be applied in real time.
One can also use Equation \eqref{eq:mainprob} for a range of values
of $\gamma$ to obtain the probability of being contained in each member of a family
of super-level sets, letting one associate probabilities with all points in state
space. 

These two estimates also provide information about both individual consensus runs
and families of consensus runs. Specifically, Equation \eqref{eq:nmu} applies to single
consensus runs and evolves the current expected decrease in $V$ based upon
$\hz(k)$, providing a rate estimate specific to that run.
On the other hand, Equation \eqref{eq:mainprob} applies to all runs
starting from a given initial condition, giving information about how often we should
expect one trajectory out of a family to be at least some distance from the point
of convergence. Furthermore, they are complementary in that Equation \eqref{eq:nmu}
provides an upper bound on the rate of decrease in $V$ and is optimistic in the
sense that it over-estimates the expected decrease in disagreement in the system.
On the other hand, Equation \eqref{eq:mainprob} over-estimates the probability
of the agents' disagreement being a certain size, and is therefore pessimistic. 

Taken together, these two convergence rate estimates enable one to probe the behavior
of any consensus problem over random graphs, regardless of network size
$n$, edge probability $p$, or initial condition $z(0)$, and provide
quantitative data on such problems while requiring only the computation of $\mu(n, p)$. 
In the next section we present numerical results that verify both bounds presented in Theorem~\ref{thm:main}. 

\begin{figure}
\centering
\includegraphics[width=3.6in]{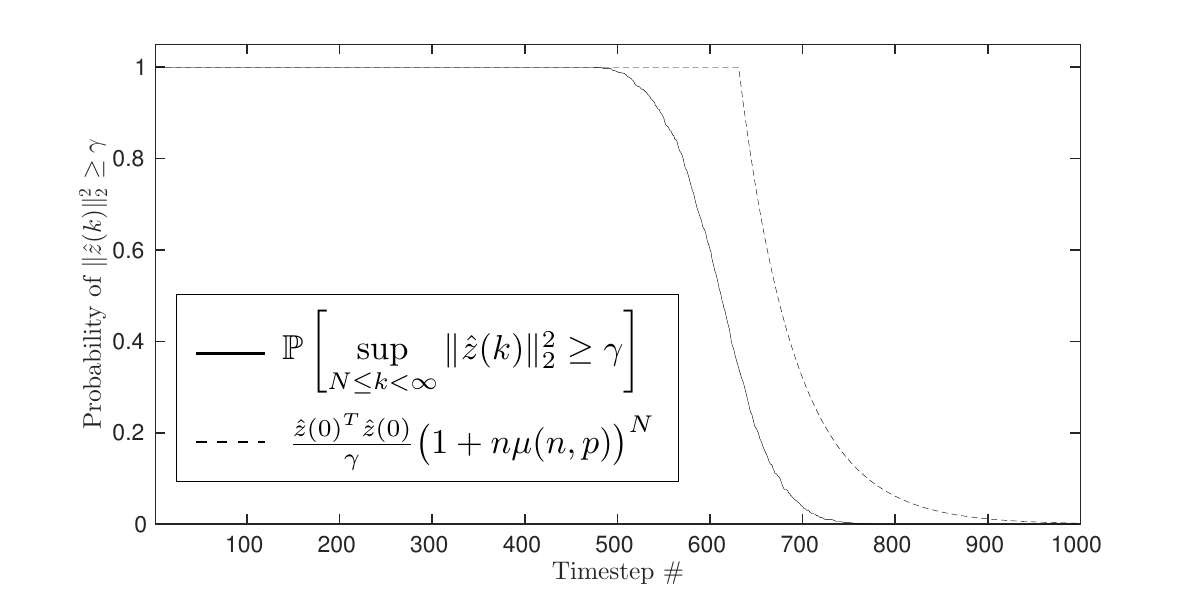}
\caption{A simulation run showing the probability of having 
$\|\hz(k)\|_2^2 \geq 3$. The empirical probability is plotted as the
solid line, while the upper bound is
shown as the dashed line. We see that the upper bound indeed holds
and, despite relying only on expected values, provides a 
close enough approximation to the empirical probability to be useful
in a variety of settings. 
}
\label{fig:prob}
\end{figure}

%


\section{Simulation Results} \label{sec:simulation}


In this section we present numerical results to support the results in Theorem~\ref{thm:main}.
We simulate consensus over random graphs and first numerically examine the expected
decrease in $V$,
and second bound the probability of being at least some distance away from the point
of agreement.

\subsection{Estimating Decreases in Disagreement}
We now present simulation results to verify the upper bound
on the expected decrease in the Lyapunov function $V(z) = \frac{1}{n}z^T\Lhat z$
as defined in Equation~\eqref{eq:Vdef}. The consensus problem we ran consists
of $n = 50$ agents and edge probability $p = 0.03$. All agents had two states
and were initialized to be evenly spaced along a circle of radius $100$ centered
on the origin. In this case, to estimate the rate of convergence using Equation~\eqref{eq:nmu},
we compute $n\mu(n, p) = -0.0561$, from which we find
\begin{equation} \label{eq:sim1}
\E\big[V\big(\hz(k+1)\big) - V\big(\hz(k)\big) \mid \hz(k)\big] \leq -0.0561\|\hz(k)\|_2^2.
\end{equation}

To validate Equation~\eqref{eq:sim1}, a single consensus run was simulated. 
At each timestep $k$, the value of $\E\big[V\big(\hz(k+1)\big) - V\big(\hz(k)\big) \mid \hz(k)\big]$ was computed numerically
by fixing $\hz(k)$, generating $1,000$ random graphs to compute
$\tilde{z}(k+1) = e^{-2\delta L(G)}\hz(k)$ for each graph $G$ generated this way, and 
then computing $V\big(\tilde{z}(k+1)\big) - V\big(\hz(k)\big) \mid \hz(k)$ for each $G$. These values were then
averaged to numerically determine $\E\big[V\big(\hz(k+1)\big) - V\big(\hz(k)\big) \mid \hz(k)\big]$ before
the algorithm proceeded to run and computed $\hz(k+1)$. 

The results of this simulation run are shown in Figure~\ref{fig:vdiff}, where 
the right-hand side of Equation~\eqref{eq:sim1} is shown as a dashed line and the
left-hand side of Equation~\eqref{eq:sim1} is shown as a solid line; though
$1,000$ timesteps of consensus were run, only the first $100$ are shown because the
two lines are indistinguishable and approximately zero beyond this point. 
Figure~\ref{fig:vdiff} shows that indeed the upper bound on the decreases in
$V$ from Theorem~\ref{thm:main} holds because the dashed line is always
above the solid line. 
Furthermore, as the iteration count
increases, the upper bound becomes more accurate, meaning that we not
only have an upper bound on the rate of decrease of $V$, but also that
this upper bound can accurately predict when decreases in $V$ go to zero,
thereby accurately predicting when consensus is achieved. 

\subsection{Bounding the Probability of Being away from $\bar{z}(0)\bone$}
We now present a consensus problem to verify the upper bound on the probability 
that the disagreement among agents will be at least a certain amount
after a fixed point in time, stated in Theorem~\ref{thm:main}. 
In particular, for $n = 10$ agents and edge probability $p = 0.01$,
$1,000$ trials were run to find experimentally the value of
$\mathbb{P}\big[\sup_{N \leq k < \infty}\|\hat{z}(k)\|_2^2 \geq \gamma\big]$
for each $N \in \{1, \ldots, 1,000\}$, where we set $\gamma = 3$. 
All trials were initialized with the agents spaced equally along 
a circle of radius $100$ whose center was at the origin.
The results of these numerical experiments are shown in Figure~\ref{fig:prob}
in the solid line, while the theoretical upper bound
$\frac{\hz(0)^T\hz(0)}{\gamma}\big(1 + n\mu(n, p)\big)^N$
is shown in Figure~\ref{fig:prob} as the dashed line. 

Figure~\ref{fig:prob} shows that the probability bound in Theorem~\ref{thm:main}
indeed holds because the dashed line is always aligned with or above
the solid line. In addition, we see that the upper bound's graph over time
stays close to that of the empirical probability, indicating that,
despite relying only on expected values, the bound on the probability
of being at least some distance from consensus provides a useful
estimate of the actual probability, enabling one to make
predictions about the magnitude of disagreement in a network over time.

\section{Conclusion}
Explicit convergence rate bounds
were presented for consensus over random graphs. A key
feature was that convergence rate estimates
were given in terms of the network size and edge probability
without making any assumptions about either. Eigenvalues
of the expected exponential of random graphs' Laplacians
were computed and used to derive approximate convergence
rate bounds. Numerical results 
confirmed that these results accurately capture the behavior of
consensus over random graphs.

\section*{Acknowledgments}
The authors would like to thank Professor 
Daniel Spielman for his helpful comments
on this work and for steering us to the relevant
literature on spectral graph theory. 

\bibliographystyle{plain}{}
\bibliography{sources}

\end{document}